\documentclass[leqno,10pt]{article}

\input{amssym}
\begin{document}
\title{{\bf Affine Classification of n-Curves}}
\date{}
\author{Mehdi Nadjafikhah \and Ali Mahdipour Sh.}
\maketitle
\begin{abstract} Classification of curves up to affine
transformation in a finite dimensional space was studied by some
different methods. In this paper, we achieve the exact formulas of
affine invariants via the equivalence problem and in the view of
Cartan's theorem and then, state a necessary and sufficient
condition for classification of n--Curves.
\end{abstract}
\noindent {\bf A.M.S. 2000 Subject Classification}: 53A15, 53A04,
53A55.\\
\noindent {\bf Key words}: Affine differential geometry, curves in
Euclidean space, differential invariants.
\section{Introduction}
This paper devoted to the study of curve invariants, in an
arbitrary finite dimensional space, under the group of special
affine transformations. This work was done before in some
different methods. Furthermore, these invariants were just pointed
by Spivak \cite{Sp}, in the method of Cartan's theorem, but they
were not determined explicitly. Now, we will exactly determine
these invariants in the view of Cartan's theorem and equivalence
problem.\par
An {\it affine transformation}, in a $n-$dimensional space, is
generated by the action of the general linear group ${\rm
GL}(n,{\Bbb R})$ and then, the translation group ${\Bbb R}^n$. If
we restrict ${\rm GL}(n,{\Bbb R})$ to special linear group ${\rm
SL}(n,{\Bbb R})$ of matrix with determinant equal to 1, we have a
{\it special} affine transformation. The group of special affine
transformations has $n^2+n-1$ parameters. This number is also, the
dimension of Lie algebra of special affine transformations Lie
group. The natural condition of differentiability is ${\cal
C}^{n+2}$.\par In next section, we state some preliminaries about
Maurer--Cartan forms, Cartan's theorem for the equivalence
problem, and a theorem about number of invariants in a space. In
section three, we obtain the invariants and then with them, we
classify the $n-$curves of the space.
%
\section{Preliminaries}
Let $G\subset{\rm GL}(n,{\Bbb R})$ be a matrix Lie group with Lie
algebra ${\goth g}$ and $P:G\rightarrow Mat(n\times n)$ be a
matrix-valued function which embeds $G$ into $Mat(n\times n)$ the
vector space of $n\times n$ matrices with real entries. Its
differential is $dP_B: T_BG\rightarrow T_{P(B)}Mat(n\times
n)\simeq Mat(n\times n)$.
\paragraph{Definition 2.1} The following form of $G$ is called {\it
Maurer-Cartan form}:
$$\omega_B=\{P(B)\}^{-1}\,.\,dP_B$$
that it is often written $\omega_B = P^{-1}\,.\,dP$.
The Maurer-Cartan form is the key to classifying maps into
homogeneous spaces of $G$, and this process needs to this theorem
(for proof refer to \cite{Iv-La}):
\paragraph{Theorem 2.2 (Cartan)}{\it Let $G$ be a matrix Lie group with Lie
algebra ${\goth g}$ and Maurer-Cartan form $\omega$. Let $M$ be a
manifold on which there exists a ${\goth g}-$valued 1-form $\phi$
satisfying $d\,\phi = -\phi\wedge\phi$. Then for any point $x\in
M$ there exist a neighborhood $U$ of $x$ and a map $f:U\rightarrow
G$ such that $f^{\ast}\,\omega = \phi$. Moreover, any two such
maps $f_1, f_2$ must satisfy $f_1 = L_B\circ f_2$ for some fixed
$B\in G$ ($L_B$ is the left action of $B$ on $G$).}
\paragraph{Corollary 2.3}{\it Given maps $f_1, f_2:M\rightarrow
G$, then $f_1^{\ast}\,\omega =f_2^{\ast}\,\omega$, that is, this
pull-back is invariant, if and only if $f_1 = L_B\circ f_2$ for
some fixed $B\in G$.} \par

\medskip The next section, is devoted to the study of the properties
of $n-$curves invariants, under the special affine transformations
group. The number of essential parameters (dimension of the Lie
algebra) is $n^2+n-1$. The natural assumption of differentiability
is ${\cal C}^{n+2}$.\par
We achieve the invariants of the $n-$curve in respect to special
affine transformations, and with theorem 2.2, two $n-$curves in
${\Bbb R}^n$ will be equivalent under special affine
transformations, if they differ with a left action introduced with
an element of ${\rm SL}(n,{\Bbb R})$ and then a translation.
%
%
\section{Classification of $n-$curves}
Let $C:[a,b]\rightarrow{\Bbb R}^n$ be a curve of class ${\cal
C}^{n+2}$ in finite dimensional space ${\Bbb R}^n$, {\it
$n-$space}, which satisfies in the condition
\begin{eqnarray}
\hspace{2cm}&\det(C',C'',\cdots,C^{(n)})\neq 0, \label{eq:3.1}
\end{eqnarray}
that, we call this curve {\it $n-$curve}. The condition
(\ref{eq:3.1}) Guarantees that $C'$, $C''$, $\cdots$, and
$C^{(n)}$ are independent, and therefore, the curve does not turn
into the lower dimensional cases. Also, we may assume that
\begin{eqnarray}
\det(C',C'',\cdots,C^{(n)})>0,
\end{eqnarray}
 to being avoid writing the absolute value in computations.\par
For the $n-$curve $C$, we define a new curve
$\alpha_C(t):[a,b]\rightarrow {\rm SL}(n,{\Bbb R})$ that is in the
following form
 \begin{eqnarray}
\alpha_C(t):=\frac{(C',C'',\cdots,C^{(n)})}{\sqrt[n]{\det(C',C'',\cdots,C^{(n)})}}\:.
 \end{eqnarray}
Obviously, it is well-defined on $[a,b]$. We can study this new
curve in respect to special affine transformations, that is the
action of affine transformations on first, second, ..., and
$n^{\textit{th}}$ differentiation of $C$. For $A$, the special
affine transformation, there is a unique representation
$A=\tau\circ B$ which $B$ is an element of ${\rm SL}(n,{\Bbb R})$
and $\tau$ is a translation in ${\Bbb R}^n$. If two $n-$curves $C$
and $\bar{C}$ be same under special affine transformations, that
is, $\bar{C}=A\circ C$, then from \cite{O'n}, we have
\begin{eqnarray}
\bar{C}'=B\circ C', \;\;\; \bar{C}''=B\circ C'',\;\;\;
...\;,\;\;\;  \bar{C}^{(n)}=B\circ C^{(n)}.
 \end{eqnarray}
We can relate the determinants of these curves as below
\begin{eqnarray}
  \det(\bar{C}',\bar{C}'',\cdots,\bar{C}^{(n)})&  = & \det(B\circ C',B\circ C'',\cdots,B\circ\bar{C}^{(n)})\nonumber \\
  & = &\det(B\circ(C',C'',\cdots,C^{(n)}))\\
 & = &\det(C',C'',\cdots,C^{(n)}). \nonumber
\end{eqnarray}
So we can conclude that $\alpha_{\bar{C}}(t)=B\circ\alpha_C(t)$
and thus $\alpha_{\bar{C}}=L_B\circ\alpha_C$ that $L_B$ is a left
translation by $B\in{\rm SL}(n,\Bbb{R})$.\par This condition is
also necessary because when $C$ and $\bar{C}$ are two curves in
$\Bbb{R}^n$ such that for an element $B\in {\rm SL}(n,\Bbb R)$, we
have $\alpha_{\bar{C}}=L_B\circ\alpha_C$, thus we can write
 \begin{eqnarray}
  \alpha_{\bar{C}}(t) &=& \det(\bar{C'},\bar{C''},\cdots,\bar{C}^{(n)})^{-1/n}(\bar{C'},\bar{C''},
  \cdots,\bar{C}^{(n)}) \nonumber \\
  &=& \det(B\circ(C',C'',\cdots,C^{(n)}))^{-1/n} B\circ(C',C'',\cdots,C^{(n)}) \\
  &=& \det(C',C'',\cdots,C^{(n)})^{-1/n} B\circ(C',C'',\cdots,C^{(n)}).\nonumber
\end{eqnarray}
Therefore, we have $\bar{C}'=B\circ C'$, and so there is a
translation $\tau$ such that $A=\tau\circ B$, and so, we have
$\bar{C}=A\circ C$ when, $A$ is a $n-$dimensional affine
transformation. Therefore, we have
\paragraph{Theorem 3.1} {\it Two $n-$curves $C$ and $\bar{C}$ in ${\Bbb R}^n$ are same under
the special affine transformations that is, $\bar{C}=A\circ C$,
which $A=\tau\circ B$ for translation $\tau$ in ${\Bbb R}^n$ and
$B\in{\rm SL}(n,{\Bbb R})$; if and only if,
$\alpha_{\bar{C}}=L_B\circ\alpha_C$, where $L_B$ is left
translation by $B$.}\\

From Cartan's theorem, a necessary and sufficient condition for
$\alpha_{\bar{C}}=L_B\circ\alpha_C$ by $B\in {\rm SL}(n,{\Bbb
R})$, is that for any left invariant 1-form $\omega^i$ on ${\rm
SL}(n,{\Bbb R})$ we have
$\alpha_{\bar{C}}^{\ast}(\omega^i)=\alpha_{C}^{\ast}(\omega^i)$,
that is equivalent with
$\alpha_{\bar{C}}^{\ast}(\omega)=\alpha_{C}^{\ast}(\omega)$, for
natural ${\goth sl}(n,{\Bbb R})$-valued 1-form
$\omega=P^{-1}\,.\,dP$, where $P$ is the Maurer--Cartan form.\par
Thereby, we must compute the $\alpha_C^{\ast}(P^{-1}.dP)$, which
is invariant under special affine transformations, that is, its
entries are invariant functions of $n-$curves. This $n\times n$
matrix form, consists of arrays that are coefficients of $dt$.

\medskip Since
$\alpha_C^{\ast}(P^{-1}\,.\,dP)=\alpha_C^{-1}\,.\,d\alpha_C$, so
for finding the invariants, it is sufficient that we calculate the
matrix $\alpha_C(t)^{-1}$. $d\alpha_C(t)$. Thus, we compute
$\alpha_C^{\ast}(P^{-1}\,.\,dP)$. We have
\begin{eqnarray}
\alpha_C^{-1}=\sqrt[n]{\det(C',C'',\cdots,C^{(n)})}\,.\,(C',C'',\cdots,C^{(n)})^{-1}.
\end{eqnarray}
We assume that $C$ is in the form $(C_1 \;\; C_2 \;\;\cdots\;\;
C_n)^T$. By differentiating of determinant, we have
\begin{eqnarray}
[\det(C',C'',\cdots,C^{(n)})]' &=&\det(C'',C'',\cdots,C^{(n)})\nonumber \\
&& + \det(C',C''',\cdots,C^{(n)}) \nonumber \\
&& \vdots\\
&& + \det(C',C'',\cdots,C^{(n-1)},C^{(n+1)})\nonumber \\
&=&\det(C',C'',\cdots,C^{(n-1)},C^{(n+1)})\,.\nonumber
\end{eqnarray}
Thus, we conclude that
\begin{eqnarray}
\alpha'_C&=&\{\det(C',C'',\cdots,C^{(n)})\}^{-1/n}\,\cdot
\left(\begin{array}{cccc}
C''_1 & C'''_1  & \cdots & C^{(n)}_1 \\
C''_ 2& C'''_ 2 &  \cdots & C^{(n)}_2\\
\vdots & \vdots       &         & \vdots \\
 C''_3 & C'''_ 3&  \cdots & C^{(n)}_n)
\end{array}\right)\nonumber \\ && \\
&& -\frac{1}{n}\,\det(C',C'',C''')\}^{-(n+1)/n}\,\cdot \left(
\begin{array}{cccc}
C'_1 & C''_1 & \cdots & C^{(n)}_1 \\
C'_ 2& C''_ 2 & \cdots & C^{(n)}_ 2\\
\vdots & \vdots       &         & \vdots \\
 C'_3 & C''_ 3 & \cdots & C^{(n)}_3
\end{array}\right).\nonumber
\end{eqnarray}
Therefore, we have the $\alpha_C^{-1} \cdot d\alpha_C$ as the
following matrix multiplying with $dt$:
\begin{eqnarray}
\left.\left(\begin{array}{ccccc}
a & 0  & \cdots & 0      & 0 \\
1 & a  & \cdots & 0      & 0 \\
\vdots & \vdots & \ddots & \vdots & \vdots \\
0 & 0 & \cdots & 1 & a \\
0 & 0 & \cdots & 0 & 1
\end{array}\right|
M\,.\,C^{(n+1)} + \left(\begin{array}{c} 0 \\ 0 \\ \vdots
\\ a \end{array}\right)\right)
\end{eqnarray}
where, the latest column, $M\,.\,C^{(n+1)} + (0,0,\cdots,a)^T$ is
multiple of $M$ by $C^{(n+1)}$ added by transpose of
$(0,0,\cdots,a)$; which $M$ is the inverse of matrix
$(C',C'',\cdots,C^{(n)})$ and also, we assumed that
\begin{eqnarray}
a=-\frac{\det(C',C'',\cdots,C^{(n-1)},C^{(n+1)})}{n\;\det(C',C'',\cdots,C^{(n)})}.
\end{eqnarray}
But with use of Crammer's law, we compute $M\,.\,C^{(n+1)}$. If
$M\,.\,C^{(n+1)} = X =(X_1,X_2,\cdots,X_n)^T$, then $M^{-1}.X =
C^{(n+1)}$. So for each $i=1,2,\cdots,n$ we conclude that
\begin{eqnarray}
X_i=\frac{\det(C',C'',\cdots,C^{(i-1)},C^{(n+1)},C^{(i+1)},\cdots,C^{(n)})}{\det(C',C'',\cdots,C^{(n)})}
\end{eqnarray}
Finally, the $\alpha_C^{-1} \cdot d\alpha_C$ is the following
multiple of $dt$:
\begin{eqnarray}
\left( \begin{array}{cccccc} a & 0 & \cdots & 0 & 0
& (-1)^{n-1}\,\frac{\det(C'',\cdots,C^{(n+1)})}{\det(C',C'',\cdots,C^{(n)})} \\
1 & a & \cdots & 0 & 0
& (-1)^{n-2}\,\frac{\det(C',C''',\cdots,C^{(n)})}{\det(C',C'',\cdots,C^{(n)})} \\
\vdots & \vdots & \ddots & \vdots  & \vdots & \vdots \\
o & 0 & \cdots & 1 & a
& -\frac{\det(C',\cdots,C^{(n-2)},C^{(n)},C^{(n+1)})}{\det(C',C'',\cdots,C^{(n)})} \\
0 & 0 & \cdots & 0 & 1 &
\frac{n-1}{n}\frac{\det(C',\cdots,C^{(n-1)},C^{(n+1)})}{\det(C',C'',\cdots,C^{(n)})}
\end{array}\right), \label{eq:3.2}
\end{eqnarray}
which the coefficient $(-1)^{i-1}$ for $i^{\textit{th}}$ entry of
the latest column, comes from the translation of $C^{(n+1)}$ to
$n^{\textit{th}}$ column of matrix
\begin{eqnarray}
(C',C'',\cdots,C^{(i-1)},C^{(n+1)},C^{(i+1)},\cdots,C^{(n)}).
\end{eqnarray}
Clearly, the trace of matrix (\ref{eq:3.2}) is zero. The entries
of $\alpha_C^{\ast}(P^{-1}\,.\,dP)$ and therefore, arrays of
matrix (\ref{eq:3.2}), are invariants of the group action.\par
Two $n-$curves $C,\bar{C}:[a,b]\rightarrow{\Bbb R}^n$ are same in
respect to special affine transformations, if we have
\begin{eqnarray}
\frac{\det(C''(t),\cdots,C^{(n+1)}(t))}{\det(C'(t),C''(t),\cdots,C^{(n)}(t))}&=&
\frac{\det(\bar{C}''(t),\cdots,\bar{C}^{(n+1)}(t))}{\det(\bar{C}'(t),\bar{C}''(t),\cdots,\bar{C}^{(n)}(t))}\nonumber \\[2mm]
\frac{\det(C'(t),C'''(t),\cdots,C^{(n+1)}(t))}{\det(C'(t),C''(t),\cdots,C^{(n)}(t))}
&=&
\frac{\det(\bar{C}'(t),\bar{C}'''(t),\cdots,\bar{C}^{(n+1)}(t))}{\det(\bar{C}'(t),\bar{C}''(t),\cdots,\bar{C}^{(n)}(t))}\nonumber \\
 & \vdots &\\
\frac{\det(C'(t),\cdots,C^{(n-1)}(t),C^{(n+1)})}{\det(C'(t),C''(t),\cdots,C^{(n)}(t))}&=&
\frac{\det(\bar{C}'(t),\cdots,\bar{C}^{(n-1)}(t),\bar{C}^{(n+1)})}{\det(\bar{C}'(t),\bar{C}''(t),\cdots,\bar{C}^{(n)}(t))}.\nonumber
\end{eqnarray}
We may use of a proper parametrization
$\gamma:[a,b]\rightarrow[0,l]$, such that the parameterized curve,
$\gamma=C\circ\sigma^{-1}$,
 satisfies in condition
\begin{eqnarray}
\det(\gamma'(s),\gamma''(s),\cdots,\gamma^{(n-1)}(s),\gamma^{(n+1)}(s))=0,
\end{eqnarray} then, the arrays
 on main diagonal of $\alpha_\gamma^{\ast}(dP\,.\,P^{-1})$ will be
 zero. But the latest determinant is differentiation of
 $\det(\gamma'(s),\gamma''(s),\cdots,\gamma^{(n)}(s))$ thus, it is
 sufficient that we suppose
\begin{eqnarray}
\det(\gamma'(s),\gamma''(s),\cdots,\gamma^{(n)}(s))=1.
\end{eqnarray}
 On the other hand, we have
 \begin{eqnarray}
   C' &=&(\gamma\circ\sigma)' = \sigma'.(\gamma'\circ\sigma) \nonumber \\
   C'' &= &
   (\sigma')^2.(\gamma''\circ\sigma)+\sigma''.(\gamma'\circ\sigma)\nonumber \\
   & \vdots&\\
   C^{(n)} &=& (\sigma')^{(n)}.(\gamma^{(n)}\circ\sigma)+ n\,\sigma^{(n-1)}\sigma'.(\gamma^{(n-1)}\circ\sigma)\nonumber \\
   & & +\frac{n(n-1)}{2}\,\sigma^{(n-2)}\sigma''.(\gamma^{(n-2)}\circ\sigma) + \cdots
   + \sigma^{(n)}.(\gamma'\circ\sigma)\nonumber
\end{eqnarray}
Therefore, $C^{(i)}$s for $1\leq i\leq n$, are some statements in
respect to $\gamma^{(j)}\circ\sigma$, $1\leq j\leq n$. We conclude
that
\begin{eqnarray}
\det(C',C'',\cdots,C^{(n)})\!\!\! &=&\!\!\!
  \det(\sigma'.(\gamma'\circ\sigma)\,,\,(\sigma')^2.(\gamma''\circ\sigma)+\sigma''.(\gamma'\circ\sigma)\,,\,\nonumber \\
  & & \cdots, (\sigma')^{(n)}.(\gamma^{(n)}\circ\sigma)+ n\,\sigma^{(n-1)}\sigma'.(\gamma^{(n-1)}\circ\sigma)\nonumber \\
  && + \frac{n(n-1)}{2}\,\sigma^{(n-2)}\sigma''.(\gamma^{(n-2)}\circ\sigma) +
  \cdots + \sigma^{(n)}.(\gamma'\circ\sigma)) \\
  & = & \det(\sigma'.(\gamma'\circ\sigma)\,,\,(\sigma')^2.(\gamma''\circ\sigma)\,,\,\cdots\,,\,(\sigma')^n.(\gamma^{(n)}\circ\sigma))\nonumber \\
  & = & \sigma'^{\frac{n(n-1)}{2}} . \det(\gamma'\circ\sigma\,,\,\gamma''\circ\sigma\,,\,\cdots\,,\,\gamma^{(n)}\circ\sigma)\nonumber \\
  & = & \sigma'^{\frac{n(n-1)}{2}}\,,\nonumber
\end{eqnarray}
The latest expression signifies  $\sigma$ therefore, we define the
{\it special affine arc length} as follows
\begin{eqnarray}
\sigma(t):=\int_a^t\;\Big\{\det(C'(u),C''(u),\cdots,C^{(n)}(u))\Big\}^{\frac{2}{n(n-1)}}\,du.
\end{eqnarray}
So, $\sigma$ is the natural parameter for $n-$curves under the
action of special affine transformations, that is, when $C$ be
parameterized with $\sigma$, then for each special affine
transformation $A$, $A\circ C$ will also be parameterized with the
same $\sigma$. Furthermore, every $n-$curve parameterized with
$\sigma$ in respect to special affine transformations, will be
introduced with the following invariants
\begin{eqnarray}
 \chi_1    &=&(-1)^{n-1}\,\det(C'',\cdots,C^{(n+1)})\nonumber \\
 \chi_2    &=&(-1)^{n-2}\,\det(C',C''',\cdots,C^{(n)})\nonumber \\
&\vdots& \\
 \chi_{n-1} &=& \det(C',\cdots,C^{(n-2)},C^{(n)},C^{(n+1)}).\nonumber
\end{eqnarray}

We call $\chi_1 ,\chi_2, \cdots ,$ and $\chi_{n-1}$ as
(respectively) the first, second, ..., and $n-1^{\textit{th}}$
{\it special affine curvatures}. In fact, we proved the following
theorem
\paragraph{Theorem 3.2}{\it
A curve of class ${\cal C}^{n+2}$ in ${\Bbb R}^n$ with condition
(3.1), up to special affine transformations has $n-1$ invariants
$\chi_1$, $\chi_2$, ..., and $\chi_{n-1}$, the first, second, ...,
and $n-1^{\textit{th}}$ affine curvatures that are defined as
formulas (3.3)}.
\paragraph{Theorem 3.3}{\it
Two $n-$curves $C, \bar{C}: [a,b]\rightarrow{\Bbb R}^n$ of class
${\cal C}^{n+2}$, that satisfy in the condition (3.1), are special
affine equivalent, if and only if, $\chi_1^C=\chi_1^{\bar{C}}$,
$\cdots$, and $\chi_{n-1}^C=\chi_{n-1}^{\bar{C}}$.}\\

\noindent{\it Proof:} Proof is completely similar to the three
dimensional case \cite{Nm}. The first side of the theorem was
proved in above descriptions. For the other side, we assume that
$C$ and $\bar{C}$ are $n-$curves of class ${\cal C}^{n+2}$ with
conditions (resp.):
\begin{eqnarray}
\label{eq:3.4}
\det(C',C'',\cdots,C^{(n)})>0,\hspace{1cm}
\det(\bar{C}',\bar{C}'',\cdots,\bar{C^{(n)}})>0,
\end{eqnarray}
with this mean that they are not $(n-1)$--curves. Also, we suppose
that they have same $\chi_1$, $\cdots$, and $\chi_{n-1}$.\par
By changing the parameter to the natural parameter ($\sigma$),
discussed above, we obtain new curves $\gamma$ and $\bar{\gamma}$
resp. that the determinants (\ref{eq:3.4}) will be equal to 1. We
prove that $\gamma$ and $\bar{\gamma}$ are special affine
equivalent, so there is a special affine transformation $A$ such
that $\bar{\gamma}=A\circ\gamma$ and then we have $\bar{C}=A\circ
C$ and proof will be completed.\par
At first, we replace the curve $\gamma$ with
$\delta:=\tau(\gamma)$ properly, in which case that $\delta$
intersects $\bar{\gamma}$, that $\tau$ is a translation defined by
translating one point of $\gamma$ to one point of $\bar{\gamma}$.
We correspond $t_0\in [a,b]$, to the intersection point of
$\delta$ and $\bar{\gamma}$ thus, $\delta(t_0)=\bar{\gamma}(t_0)$.
One can find a unique element $B$ of the general linear group
${\rm GL}(n,\Bbb{R})$, such that maps the base
$\{\delta'(t_0),\delta''(t_0),\cdots,\delta^{(n)}(t_0)\}$ of
tangent space $T_{\delta(t_0)}{\Bbb R}^3$ to the base
$\{\bar{\gamma}'(t_0),\bar{\gamma}''(t_0),\cdots,\bar{\gamma}^{(n)}(t_0)\}$
of it. So, we have $B\circ\delta'(t_0)=\bar{\gamma}'(t_0)$,
$B\circ\delta''(t_0)=\bar{\gamma}''(t_0)$, $\cdots$, and
$B\circ\delta^{(n)}(t_0)=\bar{\gamma}^{(n)}(t_0)$. $B$ also is an
element of the special linear group, ${\rm SL}(n,\Bbb{R})$, since
we have
\begin{eqnarray}
&& \hspace{-1cm}
\det(\gamma'(t_0),\gamma''(t_0),\cdots,\gamma^{(n)}(t_0))=
\nonumber \\ && \hspace{1cm}
=\det(\delta'(t_0),\delta''(t_0),\cdots,\delta^{(n)}(t_0)),
\end{eqnarray}
and
\begin{eqnarray}
&& \hspace{-1cm}
\det(\delta'(t_0),\delta''(t_0),\cdots,\delta^{(n)}(t_0))=\nonumber \\
&& \hspace{1cm} \det\Big(B \circ
(\bar{\gamma}'(t_0),\bar{\gamma}''(t_0),\cdots,\bar{\gamma}^{(n)}(t_0)
)\Big),
\end{eqnarray}
so, $\det(B)=1$. If we denote that $\eta:=B\circ\delta$ is equal
to $\bar{\gamma}$ on $[a,b]$, then by choosing $A=\tau\circ B$,
there will remind nothing for proof.\par
 For the curves $\eta$ and $\bar{\gamma}$ we have (resp.)
\begin{eqnarray}
&& (\eta',\eta'',\cdots,\eta^{(n)})'=\nonumber \\
&& \hspace{1cm} =(\eta',\eta'',\cdots,\eta^{(n)}).\left(
\begin{array}{cccccc}
0 & 0 & \cdots & 0 & 0 &\chi_1^{\eta} \\
1 & 0 & \cdots & 0 & 0 &-\chi_2^{\eta} \\
0 & 1 & \cdots & 0 & 0 &\chi_3^{\eta} \\
\vdots & \vdots & \ddots & \vdots & \vdots & \vdots\\
0 & 0 & \cdots & 1 & 0 & (-1)^{(n-2)}\chi_{n-1}^{\eta}\\
0 & 0 & \cdots & 0 & 1 & 0
\end{array}\right),
\end{eqnarray}
and
\begin{eqnarray}
&& (\bar{\gamma}',\bar{\gamma}'',\cdots,\bar{\gamma}^{(n)})'=\nonumber\\
&& \hspace{1cm}
=(\bar{\gamma}',\bar{\gamma}'',\cdots,\bar{\gamma}^{(n)}).\left(
\begin{array}{cccccc}
0 & 0 & \cdots & 0 & 0 &\chi_1^{\bar{\gamma}} \\
1 & 0 & \cdots & 0 & 0 &-\chi_2^{\bar{\gamma}} \\
0 & 1 & \cdots & 0 & 0 &\chi_3^{\bar{\gamma}} \\
\vdots & \vdots & \ddots & \vdots & \vdots & \vdots\\
0 & 0 & \cdots & 1 & 0 & (-1)^{(n-2)}\chi_{n-1}^{\bar{\gamma}}\\
0 & 0 & \cdots & 0 & 1 & 0
\end{array}\right),
\end{eqnarray}
Since, $\chi_1$, $\cdots$, and $\chi_{n-1}$, are invariants under
special affine transformations so, we have
\begin{eqnarray}
\chi_i^{\eta}=\chi_i^{\gamma}=\chi_i^{\bar{\gamma}}, \hspace{1cm}
(i=1,\cdots,n-1).
\end{eqnarray}
Therefore, we conclude that $\eta$ and $\bar{\gamma}$ are
solutions of ordinary differential equation of degree $n+1$:
$$Y^{n+1}+(-1)^{n-1}\chi_{n-1}Y^{(n)}+\cdots+\chi_2\,Y''-\chi_1\,Y'=0,$$
where, $Y$ depends to parameter $t$. Because of same initial
conditions
\begin{eqnarray}
\eta^{(i)}(t_0)=B\circ\delta(t_0)=\bar{\gamma}^{(i)}(t_0),
\end{eqnarray}
for $i=0,\cdots,n$, and the generalization of the existence and
uniqueness theorem of solutions, we have $\eta=\bar{\gamma}$ in a
neighborhood of $t_0$, that can be extended to all $[a,b]$. \hfill
$\diamondsuit$
\paragraph{Corollary 3.4} {\it The number of invariants of special
affine transformations group acting on ${\Bbb R}^n$ is $n-1$, that
is same with results provided with other methods such as
\cite{Gu}.}\par
{}\vspace{0.45cm}

Mehdi Nadjafikhah\\
{\it Department of Mathematics,\\
Iran University of Science and Technology,
Narmak-16, Tehran, Iran.\\ E-mail:} \verb"m\_nadjafikhah@iust.ac.ir"\\

Ali Mahdipour Sh.\\
{\it Department of Mathematics,\\
Iran University of Science and Technology, Narmak-16, Tehran,
Iran.\\ E-mail:} \verb"mahdi\_psh@mathdep.iust.ac.ir"

\end{document}